\documentstyle{amsppt}
\NoRunningHeads
\TagsOnLeft

\def\P{{\Cal P}}
\def\R3{ {\Bbb R}^3}

\def\S{{\Cal S}}

\def\a{\alpha}
\def\b{\beta}
\def\c{\wedge}
\def\e{\eta}
\def\ep{\epsilon}
\def\z{\zeta}

\def\l{\lambda}
\def\m{\mu}
\def\n{ \nabla}
\def\o{\omega}
\def\p{\partial}

\def\vO{\varOmega}
\def\vp{\varphi}
\def\vw{\varpi}

\topmatter
\title     A sufficient condition for a finite-time $ L_2 $ singularity of
the 3d Euler Equation  \endtitle 
\author Xinyu He \endauthor
\affil  Mathematics Institute \\
 University of Warwick \\
 Coventry CV4 7AL, UK \\
\smallskip
 Written   on $ \ $   30/08/2002   
\endaffil

\abstract
A sufficient condition is derived for a finite-time $ L_2 $ singularity 
of the 3d incompressible Euler equations,  making appropriate assumptions
on eigenvalues of the Hessian of pressure. 
Under this condition  $ \ \lim_{ t \to T_*} \sup
\left\|  \frac{ D \o} { Dt} \right\|_{L_2(\vO)} = \infty $,
where $ ~ \vO \subset \R3 $ moves with the fluid.  In particular, 
$  |  \o | $, $ | \S_{ij} | $ , and 
 $  | \P_{ij} | $ all become unbounded at one 
 point $ (x_1, T_1) $, $ T_1 $ being the first blow-up time in $ L_2 $.

\endabstract
\endtopmatter
\document

\noindent {\bf 1. Introduction} 
\bigskip

\noindent 
Consider the incompressible Euler equations in $ \R3 \times [0,\infty) $ 
$$ \frac{ \p  u }{ \p t} + u \cdot \n  u = - \n p , 
  \ \ \ \   \n \cdot u = 0,   \tag 1.1 $$ 
 where $~ u(x,t)=(u_1,u_2,u_3)$ denotes 
 the unknown velocity field, $ p $ the pressure scalar. 
 Denote the material derivative in (1.1) by
 $ D / Dt = \p / \p t + u \cdot \n $, and
 the vorticity vector by  $ \o = \n \c u $, which is governed by  
 $$  \frac{ D  \o }{ D  t} = \S ~ \o ,
      \ \ \ \  \n \cdot \o = 0 ,  \ \ \ \ \text{where} \ \ \ 
   \S_{ij} := \frac {1}{2 } \left( \frac {\p u_i }{ \p x_j }+ 
  \frac { \p u_j }{  \p x_i} \right). \tag 1.2  $$
 Defining  the Hessian of pressure $  p $  by
 $$ \P_{ij} := \frac { \p^2 p }{ \p x_i \p x_j } , \tag 1.3 $$
 the second order derivative of 
 $ \o $ is given  by (see [Ma] and [O]) 
 $$ \frac{ D^2  \o }{ D  t^2} = - \P ~ \o . \tag 1.4  $$

Combining (1.2) and (1.4), it is shown in [GGK] that
$$ \frac { D ( \o \c \S \o) } { D t} = - \o \c \P \o .$$
This means that if $ \o $ aligns with an eigenvector of $ \S $
(call this a  $ \S-\o$ alignment), then it must do so
simultaneously with an eigenvector of $ \P $
(call this a  $ \P-\o$  alignment).
See (3.1) for the converse.  
It is clear from (1.4) that only negative eigenvalues of $ \P $ cause 
$ \o $ to increase in time.  Intuitively, one expects that
singular solutions of (1.1), if they exist,
 are related to alignments of $ \P-\o $ or $ \S-\o $. 
 In this sense, the geometry matters.

The theorem of [BKM] states that the $ L_{\infty} $ norm of  $  \o  $ controls 
the smoothness of solutions of the Euler equations (1.1).  On the other hand, 
the direction of vorticity plays an important role with its evolution 
connected to the Hessian of pressure ~ $\P  ~ $ [ p. 40, C ]. 
It is further proved in [CFM] that
if the direction of $ \o $ remains regular and the velocity is bounded,
then a singularity cannot form. 

There has been evidence that alignments exist in a wide classes of fluid flows.
It is found in [Pe] that in the Euler singular region, 
the vorticity is aligned with
the eigenvector of the most positive eigenvalue of the strain $ \S $.
With vortex pairs initially aligned with $ \S $,
a blow-up model is constructed [Mo]. 
Using a set of equations for the angle variables in terms of $ \S $ and $ \P $,
[GGK] has recently analysed the data [K], indicative of
intense stretching and compression of vorticity at the singular region 
where the alignments occur (see their Fig. 2 and 3).
See also [H] for the alignments associated with Navier-Stokes turbulence.

The present paper is to study geometrical configurations of $ \P $.
We shall derive a sufficient condition in Theorem 2.1 for 
a finite-time $ L_2 $ Euler singularity,
assuming the direction of $ \o $ is parallel to an eigenvector of $ \P $ only.
Furthermore, assuming the direction of $ \o $ is parallel to 
both $ \P $ and $ \S $ in a simple way, Theorem 2.2 is obtained.
Deducing from this theorem,
we analyse the singular patterns in time and space by Corollary 2.3 and 2.4.
Apparently,  these patterns seem to be
observed in [K] and [Pe] for the turbulent enstrophy dissipation.
Finally, we discuss effectiveness of the Hessian of pressure on producing
potential $L_2 $ singularities.

\bigskip
\remark{Remark A }  To prove the theorems, we imposed some conditions on
the eigenvalues of $ \S $ and $ \P $.  Although little is known about a relation
between their eigenvalues, the conditions imposed may be justified by 
available numerical data.   Note that
the conditions already imply  possible pointwise Euler singularities.
 However, the central point of the paper is to 
demonstrate that a $ L_2 $ blowup demands stronger conditions. 
Our condition for a pointwise
singularity is not sufficient (see {\sl Remark D}).
Moreover, global constraints need to be satisfied,
for instance only fluid elements
satisfying inequality (2.10)  become  unbounded in $ L_2 (\vO) $.
To the author's knowledge, sufficient conditions for $ L_2 $ Euler blowup
have not been precisely derived before.

\endremark

\bigskip
\noindent {\bf 2. A sufficient condition} 

\bigskip
\noindent 
Let $ \vO \subset \R3 $ be a smooth material volume carried by the fluid. 
Let $ ~ \o(x,t) $ be a sufficiently smooth solution of (1.1) for which
we set 
$$ \vw :=   \| \o(t) \|_{L_2(\vO)}^2 , \ \ \  
\vw(t) \not= 0 \  \ \forall ~ t \geq 0 ,
\ \ \  \text{and} \ \ \ \vp_1 (t)  := \frac{1} {2 \vw } .\tag 2.1 $$
 
\medskip
\remark{Remark B }  One could also set
$$ \vp_n (t) := \frac{1}{2 ~ [ ~ \vw ~ ]^{\frac{1}{n}}}, \ \ n \in {\Bbb N} .$$
This would slightly improve an estimate for the constant $ c_0 $ in 
Theorem 2.2 below
(smaller $ c_0 $ for $ n > 1 $). 
 However for clarity, we take $ n = 1 $ as in (2.1). 
\endremark

\medskip
Define a smooth function
$$ v(t) := - ~ \vp_1'  \tag 2.2 $$  so that 
$$ v(t) = \frac {1}{\vw^2} \int_{\vO}  \ \o  \cdot  \frac{ D \o} { Dt} 
 ~ d x, \tag 2.3 $$ and
$$ v'(t) = \frac {1}{\vw^3} \left\{ 
 \left( \int_{\vO} ~ [ ~ \left| \frac{ D \o} { Dt} \right|^2
  + \o  \cdot  \frac{ D^2 \o} { Dt^2} ~ ] ~ dx \right) \vw
 \  - 4 \left ( \int_{\vO} \o  \cdot  \frac{ D \o} { Dt}  ~ d x \right)^2 
 \right\}. \tag 2.4 $$ 

\medskip
Concerning the above equations, an easy estimate is

\bigskip
\proclaim{Lemma 2.0} 
Let $ v , v' $ be as in (2.3) and (2.4).  Then for 
$ ~ t \in  [0,\infty) $ 
$$ v(t) ~ \vw^{3/2}(t) ~ \leq ~ \left\|  \frac{ D \o} { Dt} \right\|_{L_2(\vO)} ,
 \tag 2.5 $$ and
$$ v'(t) \vw^2(t) ~ \geq ~   
  \int_{\vO} ~  \o  \cdot  \frac{ D^2 \o} { Dt^2}  ~ dx
  \  - c_1  \int_{\vO} ~ \left| \frac{ D \o} { Dt} \right|^2  ~ d x  , 
  \ \ \ c_1 = 3.  \tag 2.6 $$ 
 
\endproclaim

\bigskip
\demo{Proof}  
By Cauchy-Schwarz's inequality, we get for the integral in (2.3):  
$$ \int_{\vO}  \ \o  \cdot  \frac{ D \o} { Dt} ~ dx ~ \leq ~ 
\| \o \|_{L_2(\vO)} ~  \left\|  \frac{ D \o} { Dt} \right\|_{L_2(\vO)} . $$ 
But $ \vw =   \| \o \|_{L_2(\vO)}^2 $, so this leads to (2.5).
Using this relation again for the last term in (2.4) yields (2.6).
$ \hfill{ \qed} $

\enddemo

\medskip
\remark{Remark C }
Inequality (2.6) involves both (1.2) and (1.4), therefore it will be used
to investigate various links between $ \S $ and $ \P $ for solutions of (1.1).
\endremark

\medskip
No rigorous estimate is known about the two terms on
the right hand-side of (2.6), and certain assumptions will be made
on geometrical arrangements of $ \S $ and $ \P $.  
First, we consider a case
when there is only $ \P-\o $ alignment.  This arrangement is shown by numerical 
data [O], which suggests the configuration to be a generic property of 
Euler flows.  A sufficient condition can now be given.

\bigskip
\proclaim{Theorem 2.1} 
Let  $ ~ \P \o = - \l ~ \o $ in (1.4)
$  \  \forall  ~ x \in \vO $ and $ t \geq 0 $,
where  $ \l > 0 $.  Assume that at some  $ \ t_0 > 0 $,
$ ~ \l > ~  3 ~ \m_m^2 $ on $ ~ \vO \times [t_0,\infty) $, where 
$ \m_m = \max \{ | \m_1 |,| \m_2 |,| \m_3 | \} $, 
$ \m_i $ being eigenvalues of the matrix $ \S $.
Then there exists a finite time $  T_0 > t_0 $ (depending only on $ ~ \vw_0 $ and $ ~ v_0$)
and  $ ~  T_*  \in  (t_0,  T_0) ~ $, such that   
$$ \lim_{ t \uparrow T_*}  \sup
\left\|  \frac{ D \o} { Dt} \right\|_{L_2(\vO)} = \infty .$$

\endproclaim

\bigskip
\demo{Proof} 
By Lemma 2.0, clearly
$$ v' \vw^2 ~ \geq ~   
  \int_{\vO}  ~ \l(x,t) ~ | \o |^2  ~ dx  \ -  \  
  3 \int_{\vO}  ~   | \S \o |^2  ~ dx. $$ 
  Setting $ \m_m = \max \{ | \m_i | \} $ gives 
  $$ v' \vw^2 ~ \geq ~   
  \int_{\vO} ~  \left[ ~ \l(x,t) - 3 \m_m^2(x,t) ~ \right]  ~ | \o |^2  ~ dx. $$ 
It then follows from the assumption and (2.5)
$$ v'(t)  ~ \geq ~  c ~ \vw(t) ~ v^2(t)  ,  
  \ \ \ t \in [t_0, \infty), \ \ \  c \in (0,1].  $$ 
This implies  $ \vp_1' < 0 $ in (2.2) after $ t_0  $, in turn  
$ \vw (t) ~ \ge ~ \vw_0 = \vw(t_0) $.  Hence  
$$ v'  ~ \geq ~  c ~ \vw_0 ~ v^2 ,  \ \ \ v_0 = v(t_0) > 0.  $$ 
One finds that for $ \  t_0 ~ \leq ~ t  ~ < ~ T_0 $,
setting  $ ~ A = 1/(c ~ \vw_0 ) $, 
$$ v(t) ~ \geq ~ \frac { A} { T_0 - t }, 
\ \ \ \ T_0  = t_0 + 1 / (c ~ \vw_0 ~ v_0) . $$
We see that $ ~ t_0 < T_0 < K $.
According to (2.5), in which note $ ~ \vw (t) ~ \ge ~ \vw_0  $, 
$$ \left\|  \frac{ D \o (t) } { Dt} \right\|_{L_2(\vO)} ~ \geq ~ \ 
  \frac { B} { T_0 - t}, \ \ \ \ 
 B =  \vw_0^{1/2} / c  . $$
 This establishes the assertion. $ \hfill{ \qed} $  

\enddemo

\medskip
The basic idea of Theorem 2.1 is that if $ \l $ is larger than $ \m_m $ 
for a certain length of time, then a $L_2 $ singularity forms. 
The critical time  $ T_0 $ is determined by initial $ \vw_0 $ (the enstrophy
at $ t_0$ ) and $ v_0 $ (the rate change of enstrophy); 
higher is the initial enstrophy, shorter is the critical time.

To be precise how large $ \l $ needs to be, next
we  examine a special case of the above theorem:  both $ \P - \o $ and $ \S - \o $
configurations hold.  Such flow geometry is often observed in numerical simulations,
for example [O], [GGK].  Making a assumption
on the eigenvalues of $ \S $ and $ \P $,  we have

\bigskip
\proclaim{Theorem 2.2 } 
Let  $ ~ \P \o = - \l ~ \o $ in (1.4) and $ \S \o =  \m ~ \o $ in (1.2)
$  \  \forall  ~ x \in \vO $ and $ t \geq 0 $,
where  $ \l, ~ \m > 0 $. 
 Assume that at some  $ \ t_0 > 0 $,
$ ~ \l =  c_0 \m^2 $ on $ ~ \vO \times [t_0,\infty) $ with
some constant $ ~ c_0 > 3 $.
Then there exists a finite time $ ~ T_0 > t_0 $ and 
 $ ~  T_*  \in  (t_0,  T_0) $, such that   
$$ \lim_{ t \uparrow T_*}  \sup
\left\|  \frac{ D \o} { Dt} \right\|_{L_2(\vO)} = \infty .$$

\endproclaim 

\bigskip
\demo{Proof} 
The proof is similar to that of Theorem 2.1.
Here for $ T_0 $, we have   
$$  T_0  = t_0 + 1 / (c ~ \vw_0 ~ v_0) , \ \   c = c_0 - 3 > 0  . \tag 2.7 $$
 $ \hfill{ \qed} $  

\enddemo

\medskip
\remark{Remark D } 
 When both $ \P-\o$ and $\S-\o$ alignments hold, there
 may exist many functional relations between their eigenvalues,
 $ ~ \l ~ = f(\m ) $. 
The hypothesis in the theorem, $ ~ \l =  c_0 ~ \m^2 $ with $ c_0 \in (3, 3 + \ep) $,
 is a requirement for the $ L_2 $ blowup [but note not every fluid 
 element satisfying the relation can blowup, see (2.10) below]. This
 requirement already implies pointwise singular solutions.  For 
 such singularities, a similar relation is  
$ ~ \l =  c_p ~ \m^2 $  with $ c_p \in (1, 1 + \ep) $
 (see the proof of Corollary 2.3). 
 Notice that $ ~ c_p < c_0 ~ $ for $ ~ \ep \in (0,1)$.

\endremark
\medskip

This case is the simplest to analyse structures of the $ L_2 $ blowup.
To do so we will further assume  
that $  \m   $ is the only positive eigenvalue of $ ~ \S $, as suggested
by an analysis [p. 309, Pe].  Thus 
the very first blow-up time in $ L_2 $ is identified by

\bigskip
\proclaim{Corollary 2.3 (Temporal interval)}  
Suppose in Theorem 2.2 that
$   \m   $ is the only positive eigenvalue of $ ~ \S $.
Then there exists a smallest time $ ~ T_1  \in  (t_0,  T_0) ~ $ such that
$$  \lim_{ t \uparrow T_1} \sup |  \o |_{L_{\infty}} = \infty , \ \ \ \ 
  \lim_{t \uparrow T_1} \sup | \S_{ij} |_{L_{\infty}} = \infty , \ \ \ \ 
  \text{and} \ \ \ \  \lim_{ t \uparrow T_1} \sup | \P_{ij} |_{L_{\infty}} = \infty  . $$
In fact, $ [ T_1, T_0 ) = \{ ~ t  ~ | ~ T_1 \leq t  < T_0 \} $ is the interval of blow-up. 

\endproclaim

\bigskip
\demo{Proof}
Let $ \vO_0 = \overline \vO(t_0) $ and 
 $ \m^0(x) = \m(x, t_0 ) $ for $ x \in \vO_0 $.
Consider a fluid element located at $ \a \in \vO_0 $. 
Differentiating $ D \o / Dt = \m \o $ and using (1.4), 
one obtains by following the element:
$   \m'(t)  = \l - \m^2 $.  Inserting $ ~ \l =  c_0 \m^2 $ ~ gives ~
$  \m' = (c_0 -1 ) \m^2  $. 
 This equation admits a solution which ceases to be regular at a finite-time
$$ \m(t; \a) ~  =  \frac {(c_0 -1)^{-1}} { T_* - t}, \ \ \ \ 
 ~  T_* =  t_0 + 1/ [ (c_0 -1 ) ~ \m^0(\a) ] . \tag 2.8 $$ 
Note  $ ~   \inf \m^0(\vO_0)  ~ \leq ~ \m^0(\a) 
  ~ \leq ~ \sup \m^0(\vO_0)  \ \ \forall  \ \a \in \vO_0 $.  Define
$$ ~  T_1 : = \inf_{ \a \in \vO_0 } \ T_*(\a) ~ =  t_0 + 1 / [ (c_0 -1 ) ~ \m_1^0 ], \ \ \ \ 
 \m_1^0 =  \sup \  \m^0(\vO_0) . \tag 2.9 $$

We claim $ ~ T_1 < T_0 $ as defined in (2.7).   
Computing $ ~ c~ \vw_0 ~ v_0 $  in $ T_0 $ by use of
the Second Mean-Value Theorem for Integrals in (2.3), we get
$ ~ c~ \vw_0 ~ v_0 = (c_0 -3) \m^0( \b )$ for some $ \b \in \vO_0$.
The fact  $  (c_0 -1 ) ~ \m_1^0  ~  > (c_0 -3) \m^0( \b ) 
 \  \forall  \ \b \in \vO_0 $ suffices for the claim. 
Consequently, $ T_1$ is the first time in the blow-up interval $ [ T_1, T_0 )$, 
in which corresponding $ \m^0(\a) $ necessarily satisfy
$$   \m^0(\a) ~ \geq ~ \m^0( \b_* ) ~ (c_0 -3) / (c_0 -1 ), 
   \ \  \b_* \in \vO_0 .   \tag 2.10 $$
 
We now ask what functions are singular at $ T_1 $ ?
Since both matrices  $ \S $ and $ \P $ are symmetric, we have only to 
consider their eigenvalues. 
Let $ \m_a $ and $ \m_b $ be the two other eigenvalues of $ \S $
whose eigenvectors are not aligned with the vorticity vector. 
By the incompressibility condition, 
$ \m >  \max  ~ \{ | \m_a |,  | \m_b | \} $
as it is the only positive eigenvalue.  Thus it is obvious from (2.8) 
and (2.9) that $ ~ |  \S_{ij} |_{L_{\infty}} $ is unbounded at $ T_1 $. 
  This means, by the theorems of [BKM] and [Po],
 that $ |  \o |_{L_{\infty}} $ also fails to be smooth at the same time. 
  Finally we turn to the Hessian of pressure.   
  Let $ \l_{\z} $ and $ \l_{\e} $ be the two other eigenvalues of $ \P $ 
  while $ ~ - \l ~ $ is the negative eigenvalue associated with the eigenvector aligned to $ \o $.
 Note that $ \l_{\z} $ or $ \l_{\e} $ cannot blow up at any time earlier than $ T_1 $, 
because if this happened,  it can be shown by (1.2) and (1.4)  that 
$ |  \o |_{L_{\infty}} $ would have blown up at
a time earlier than $ T_1 $, contradicting (2.9).  Now given $ ~ \delta > 0 
  ~  \forall  ~ t \in (T_1 - \delta, ~ T_1 ) $, either  
 (a) $ ~  \sup_{ x \in \vO}  ~ \l ~  \geq ~ ~  \max  ~ \{ | \l_{\z} | ,  | \l_{\e} | \}  $,
 or (b) $ ~  \sup_{ x \in \vO}  ~ \l ~  <  ~ \max  ~ \{ | \l_{\z} | ,  | \l_{\e} | \} $.
 We know that $ \lim_{t \uparrow T_1} \sup | \S_{ij} |_{L_{\infty}} = \infty $,
  which is equivalent to 
 $  \lim_{ t \uparrow T_1} ~ \sup_{ x \in \vO}  ~ \l ~  =  \infty $ by 
  the alignment relation $ ~ \l =  c_0 \m^2  $.      Thus inequality (a)
 is left as the only choice.
 Evidently  $  \lim_{ t \uparrow T_1} \sup | \P_{ij} |_{L_{\infty}} = \infty $.
The proof is complete.
$ \hfill{ \qed} $  
  
\enddemo

\medskip
 It is natural to wonder what would be the singular set in space. 
 In this direction we deduce

\bigskip
\proclaim{Corollary 2.4 (Spatial set)}
Let  $ ~ x_1 \in \vO $ be the space point where
$ |  \S_{ij} |_{L_{\infty}} = \infty  $ as $ t \to T_1 $.
Then $ |  \o |_{L_{\infty}} $  and $ |  \P_{ij} |_{L_{\infty}} $ also
blow up at $( x_1, T_1) $.

\endproclaim

\bigskip
\demo{Proof} 
Without loss of generality, let us assume that 
at time $ t_0 $, there is only one fluid element 
having  $ ~  \m_1^0 =  \sup \  \m^0(\vO_0) $. 
Suppose $ ~ |  \o |_{L_{\infty}} $  blows up at $ (y, T_1) $, $ ~ y \not= x_1 $,
 however this is impossible.
At the time $ T_1 $,  $ ~ y $ is a location arrived  by a fluid 
element with initial $ \m^0(y) \not=  \m_1^0$, which is not singular at that time.
We then conclude   $ \  y = x_1 $.  
 To find the singular location of $ | \P_{ij} |_{L_{\infty}} $,
 we recall in  Corollary 2.3 that 
 $ ~ \sup_{ x \in \vO} ~ \l ~  \geq ~ ~  \max  ~ \{ | \l_{\z} | ,  | \l_{\e} | \}  $ 
 for $ ~ t \in (T_1 - \delta, ~ T_1 ) $.
If  $ ~  \sup_{ x \in \vO}  ~ \l ~  > ~ ~  \max  ~ \{ | \l_{\z} | ,  | \l_{\e} | \}  $,
 then it is unbounded at $( x_1, T_1) $ by the alignment relation. 
 If  $ ~  \sup_{ x \in \vO}  ~ \l ~  =  ~ ~  \max  ~ \{ | \l_{\z} | ,  | \l_{\e} | \}  $,
 this means both $ ~  \sup_{ x \in \vO}  ~ \l ~ $ and 
 $ \max  ~ \{ | \l_{\z} | ,  | \l_{\e} | \}  $ blow up at $ T_1 $.
 Having stated $ ~  \sup_{ x \in \vO}  ~ \l ~ $ is singular at $( x_1, T_1) $,
 let us suppose $ \max  ~ \{ | \l_{\z} | ,  | \l_{\e} | \}  $ 
 be singular at $ (z, T_1), z \not= x_1 $.
 A similar argument to the above for $ |  \o |_{L_{\infty}} $ shows
  we must have $ ~ z = x_1 $.
$ \hfill{ \qed} $  
\enddemo

\medskip
We make a few observations about the above results.
 (i) Geometrical arrangements can limit the set of singularities.
 In the case of the double alignments, we have shown that 
 $ |  \o |$, $ |  \S_{ij} | $, and
 $ |  \P_{ij} |$ all blowup at one point 
 $( x_1, T_1) $. 
(ii)  The $L_2 $ singularity condition is stronger,  namely the integral relation 
(2.6) has to be satisfied as a constraint. 
In this instance, although in (2.8) any fluid element could locally blow up at $  T_* $, 
only those satisfying the inequality (2.10) can actually make up the $ L_2 $ singularity.
(iii) Taking the divergence of (1.1) results in
 $ ~ | \o |^2 /2 - \S^2 = \P_{ii} =  \l_{\z} +  \l_{\e}  - \l $. 
 From Corollary 2.4, we see that in any neighborhood of $(x_1, T_1)$, 
 the above equation has an indefinite sign of $ \infty - \infty $.

\bigskip
\bigskip
\noindent {\bf 3. Necessity for $ L_2 $ blow-up} 

\bigskip
 On the right hand-side of (2.6), if the 
first integral is persistently greater than the second,
then a singularity could result.
In our above theorems, we only used the geometric conditions on the integrands, which is
more restrictive than the integral requirement.  However
in general cases when there is not any coherent configuration,
it seems hard to proceed.
In what follows, we shall discuss 
solutions of (1.1) having some coherence in the Hessian of pressure.

To simplify the discussion, let  
$ \S $ and $ \P $ be diagonalised on $ ~ \vO \times [0,\infty) $
with respect to the principal axes.  
Since (2.6) is invariant under the coordinate transformations,
we can write referring to these axes 
$$ v' \vw^2 ~ \geq ~   - \int_{\vO} ~ 
 [ \l_{\z}  \o_{\z}^2  + \l_{\e}  \o_{\e}^2  + \l_{\xi} \o_{\xi}^2 ] ~ dx
  \  - 3  \int_{\vO} ~ [ \m_a^2  \o_a^2 + \m_b^2  \o_b^2 + \m_c^2  \o_c^2 ] ~ d x  ,  $$ 
where $ \z, \e , $ and $ \xi $ denote the principal axes of $ \P $,
 $a, b $ and $ c $ the principal axes of $ \S $, respectively.
It appears that a $ \P-\o$ alignment with a negative eigenvalue would be an 
effective way for attaining the requirement, for the following reason. 

As shown in the Introduction, when a $ \P-\o$ alignment occurs, we have 
$$ - \o \c \P \o \equiv  0 ~  \implies ~ \o \c \S \o ~ = \ \text{constant}. \tag 3.1$$
Let us write out three components of the invariant $ ( ~ \o \c \S \o ~ ) $:
$$ \align 
          \o_c \o_b (\m_c - \m_b) &= c_1, \\
          \o_a \o_c (\m_a - \m_c) &= c_2, \\
          \o_b \o_a (\m_b - \m_a) &= c_3  .  \tag3.2  \endalign $$ 
A key point here is that from the instant $ t_0$ at which $\P-\o$ occurs for 
some fluid elements, the constants in (3.2) are
fixed in time following the same elements. 
The configuration of a vortex tube would give 
an interesting example of (3.2).  
Suppose at $ t_0$, the fluid elements have
$ ~  \m_a > 0 $, and $ ~ \m_b ,  \m_c < 0 $ with $ \m_b = \m_c$. 
This leads to initially, $c_1 =0$, $c_2 > 0 $, and $c_3 < 0 $.
 We obtain in (3.2)
 $  \o_a = c_2/ \o_c( \m_a + | \m_c | )  $. 
  In this formula: (i) $ c_2 > 0 $ is fixed; 
 (ii) it is not clear how $ ( \m_a + | \m_c | ) $ changes in time 
 (Theorem 2.2 is not applicable); 
 (iii) $ \o_c $ decreases according to (1.2),
 since $ \m_c $ remains negative to keep $ c_1 = 0 $, due to 
 the incompressibility.  
 So there is a tendency for $ \o_a $ to increase in time, keeping
 the vortex-tube state alive, and such a state will be
 strengthened if there are some symmetries existing 
 in the flow at $ t_0$.  This (extreme) example illustrates that 
a $\P-\o$ alignment ``freezes " the initial straining states by
 (3.1), and
if the initial configuration favours vortex stretching, then these vortex lines would 
have to be stretched indefinitely. 
This suggests that the Hessian of pressure alone could
 possibly produce a $ L_2 $ singularity. 
 
The Euler equation is rich in its 
geometrical structures (see a recent paper [G]).
One further speculates whether the geometry of  $\P-\o$ or $ \S-\o$ 
is a necessary condition for solutions of (1.1) to develop 
finite-time  singularities.  Note a $ \S-\o$ alignment 
automatically implies a $ \P-\o$ alignment, but the converse is not true.
Reflecting that the alignment enforces growth of $ \o $ (cf. [Ma], p. 192), and
in view of analytical and numerical works on the subject, we may loosely make a

\bigskip
\definition {Conjecture} 
 Let $ \vO \subset \R3 $. 
Suppose (1.1) has  a $ ~ L_2(\vO) ~ $ singularity at $ ~ T_* ~ <  + \infty $.
Then $ \o $,  $ \S $, and $ ~  \P $ blow up at 
the same space point $ x_*  \in  \vO $
 $ \iff $  there exists a $ \S-\o$ alignment.

\enddefinition

\vfill\eject

\vskip  2cm
\centerline {\bf References}

\medskip\noindent 
[BKM] J.T. Beale, T. Kato and A.J. Majda, Remarks on the breakdown
of smooth solutions for the 3d Euler equations, Comm. Math. Phys.
{\bf 94} (1984), 61-66.

\medskip\noindent 
[C] P. Constantin, Geometric and analytical studies in turbulence,
in {\sl Appl. Math. Sci.} {\bf 100}, Springer-Verlag 1994.

\medskip\noindent 
[CFM] P. Constantin, C. Fefferman and A.J. Majda,
Geometric constraints on potentially singular solutions for the 3d 
Euler equations, Comm. PDE {\bf 21} (1996), 559 - 571.

\medskip\noindent
[G] J.D. Gibbon, A quaternionic structure in the 3d Euler and ideal
MHD equations, to appear in Physica D (2002).

\medskip\noindent
[GGK] J.D. Gibbon, B. Galanti and R.M. Kerr, Stretching and
compression of vorticity in the 3d Euler equations, 
in {\it Turbulence Structure and Vortex Dynamics}
 (ed. J.C.R. Hunt
\& J.C. Vassilicos), Cambridge University Press 2000.

\medskip\noindent 
[H] X. He, An invariant for the 3d Euler equations, 
 Applied Mathematics Letters {\bf 12} (1999), 55-58. 

\medskip\noindent 
[K] R.M. Kerr, Evidence for a singularity of the 3d, incompressible 
 Euler equations, Phys. Fluids A {\bf 5} (1993), 1725-1746. 

\medskip\noindent 
[Ma] A.J. Majda, Vorticity and the mathematical theory of incompressible
fluid flow, Comm. Pure Appl. Math. {\bf 39} (1986), S187-S220.

\medskip\noindent 
[Mo] H.K. Moffatt, The interaction of skewed vortex pairs:
a model for blow-up of the Navier-Stokes equations, 
J. Fluid Mech. {\bf 409} (2000), 51-68.

\medskip\noindent 
[O] K. Ohkitani, Eigenvalue problems in 3d Euler flows,
Phys. Fluids {\bf 5} (1993), 2570-2572. 

\medskip\noindent 
[Pe] R.B. Pelz, Symmetry and the hydrodynamic blow-up problem,
J. Fluid Mech. {\bf 444} (2001), 299-320 [see also: Discrete groups,
symmetric flows and hydrodynamic blowup,
in {\it Tubes, Sheets and Singularities in Fluid Dynamics}
 (ed. K. Bajer), IUTAM Symposium Series, Kluwer 2002].

\medskip\noindent 
[Po] G. Ponce, Remarks on a paper by  J.T. Beale, T. Kato, and A.J. Majda, Comm. Math. Phys.
{\bf 98} (1985), 349-353.

\enddocument
\vfill\eject
\bye